\documentclass{amsart}
\usepackage{amssymb,amscd,graphicx}

\font\nt=cmr7
\def\note#1
{\marginpar
{\nt $\leftarrow$
\par
\hfuzz=20pt \hbadness=9000 \hyphenpenalty=-100 \exhyphenpenalty=-100
\pretolerance=-1 \tolerance=9999 \doublehyphendemerits=-100000
\finalhyphendemerits=-100000 \baselineskip=6pt
#1}\hfuzz=1pt}

\def\note#1{}

\newtheorem{thm}{Theorem}[section]

\newtheorem{prop}[thm]{Proposition}

\theoremstyle{definition}

\theoremstyle{remark}
\newtheorem{rem}[thm]{Remark}
\newtheorem{ex}[thm]{Example}
\newtheorem{exs}[thm]{Exercise}

\newcommand{\C}{\mathbb C}
\newcommand{\N}{\mathbb N}
\newcommand{\Q}{\mathbb Q}
\newcommand{\R}{\mathbb R}
\newcommand{\T}{\mathbb T}
\newcommand{\Z}{\mathbb Z}
\newcommand{\no}{\noindent}
\newcommand{\dd}{\partial}
\newcommand{\sminus}{\smallsetminus}

\newcommand{\bJ}{\widetilde{\mathfrak{J}}}
\newcommand{\J}{\mathfrak{J}}
\def\<{\langle}
\def\>{\rangle}
\newcommand{\bx}{\bar{x}}
\newcommand{\bC}{\bar{C}}
\newcommand{\h}{\hbar}
\newcommand{\dvol}{dvol}
\newcommand{\cA}{\mathcal{A}}
\newcommand{\cC}{\mathcal{C}}
\newcommand{\cD}{\mathcal{D}}
\newcommand{\cF}{\mathcal{F}}
\newcommand{\cG}{\mathcal{G}}

\newcommand{\cP}{\mathcal{P}}

\newcommand{\cW}{\mathcal{W}}

\newcommand{\fG}{\mathfrak{g}}
\newcommand{\tX}{\widetilde{X}}
\newcommand{\Ga}{\alpha}
\newcommand{\Gb}{\beta}

\newcommand{\Gd}{\delta}
\newcommand{\GD}{\Delta}

\newcommand{\Ge}{\varepsilon}
\newcommand{\Gg}{\gamma}
\newcommand{\GG}{\Gamma}

\newcommand{\Gl}{\lambda}
\newcommand{\GL}{\Lambda}
\newcommand{\Go}{\omega}
\newcommand{\GO}{\Omega}
\newcommand{\Gp}{\phi}
\newcommand{\Gs}{\sigma}
\newcommand{\GS}{\Sigma}
\newcommand{\Gt}{\theta}
\newcommand{\GT}{\Theta}

\newcommand{\Aut}{\operatorname{Aut}}

\newcommand{\ev}{\operatorname{ev}}
\newcommand{\image}{\operatorname{image}}
\newcommand{\hol}{\operatorname{hol}}

\newcommand{\lk}{\operatorname{lk}}

\newcommand{\Sym}{\operatorname{Sym}}
\newcommand{\Tr}{\operatorname{Tr}}
\begin{document}
\title{{\bf Feynman diagrams for pedestrians and mathematicians}}
\author{Michael Polyak}
\address{Department of Mathematics, The Technion, 32000 Haifa, Israel}

\email{polyak@math.technion.ac.il}

\dedicatory{Dedicated to Dennis Sullivan on the occasion
of his 60th birthday}

\thanks{Partially supported by the ISF grant 86/01 and the
Loewengart research fund}



\maketitle


\section{Introduction}

\subsection{About these lecture notes}
For centuries physics was a potent source providing mathematics
with interesting ideas and problems. In the last decades something
new started to happen: physicists started to provide
mathematicians also with technical tools, methods, and solutions.
This process seem to be especially strong in geometry and
low-dimensional topology. It is enough to mention the mirror
conjecture, Seiberg-Witten invariants, quantum knot invariants,
etc.

Mathematicians, however, {\em en masse} failed to learn modern
physics. There seem to be two main obstructions. Firstly, there
are few textbooks in modern physics written in terms accessible
for mathematicians. Mathematicians and physicists speak two
different languages, and a good ``physical-mathematical dictionary"
is missing\footnote{With a notable exception of \cite{string},
which is somewhat heavy.}. Thus, to learn something from a
physical textbook, a mathematician should start from a hard and
time-consuming process of learning the physical jargon.

Secondly, mathematicians consider (and often rightly so) many
physical methods and results to be non-rigorous and do not
consider them seriously. In particular, path integrals still
remain quite problematic from a mathematical point of view (due to
some usually unclear measure aspects), so mathematicians are
reluctant to accept any results obtained by using path integrals.
Yet, this technique may be put to good use, if at least as a
tool to guess an answer to a mathematical problem.

In these notes I will focus on perturbative expansions of path
integrals near a critical point of the action. This can be done by
a standard physical technique of Feynman diagrams expansion, which
is a useful book-keeping device for keeping track of all terms in
such perturbative series. I will give a rigorous mathematical
treatment of this technique in a finite dimensional case (when it
actually belongs more to a course of multivariable calculus than to
physics), and then use a simple ``dictionary" to translate these
results to a general infinite dimensional case.

As a result, we will obtain a recipe how to write Feynman diagram
expansions for various physical theories. While in general an
input of such a recipe includes path integrals, and thus is not
well-defined mathematically, it may be used purely formally for
producing Feynman diagram series with certain expected properties.
A usual trick is then to ``sweep under the carpet" all references
to the underlying physical theory, keeping only the resulting
series. Their expected properties often can be proved rigorously,
directly from their definition.

I will illustrate these ideas on the interesting example of the
Chern-Simons theory, which leads to universal finite type
invariants of knots and 3-manifolds.

{\bf A word of caution}: during the whole treatment I will brush
aside all questions of measures, convergence, and such; see the
discussion in Section \ref{sub:divergence}.

\subsection{Basics of classical and quantum field theories}
The remaining part of this section is a brief sketch ---  on the
physical level of rigor ---  of some basic notions and physical
jargon used in the quantum field theory (QFT). Its purpose is
to give a basic mathematical dictionary of QFT's and a motivation
for our consideration of Gaussian-type integrals in this note.
An impatient reader may skip it without much harm and pass directly
to Section \ref{sec:findim}.
Good introductions to field theories can be found e.g. in \cite{DI},
\cite{Pol}; mathematical overview can be found in \cite{string};
various topological aspects of QFT are well-presented in \cite{Sc}.
Very roughly, by a {\em field theory} one usually means the following.

Given a {\em space-time} manifold $X$, one considers a space $\cF$
of {\em fields}, which are functions of some kind on $X$ (or, more
generally, sections of bundles on $X$).
A {\em Lagrangian} $L:\cF\to\R$ on $\cF$ gives rise to the
{\em action} functional $S:\cF\to\R$ defined by
$$
   S(\Gp)=\int_X L(\Gp)dx.
$$

In classical field theory one studies critical points of the
action $S$ (``classical trajectories of particles"). These fields
can be found from the variation principle $\Gd S=0$, which is simply
an infinite-dimensional version of a standard method for finding
the critical points of a smooth function $f:\R\to\R$ by solving
$\frac{d f}{d x}=0$.

In the quantum field theory one considers instead a {\em partition
function} given by a path integral
\begin{equation}
\label{eq:minkowski}
                 Z=\int_{\cF}e^{ikS(\Gp)}\cD\Gp.
\end{equation}
over the space of fields, for a constant $k\in\R$ and some formal
measure $\cD\Gp$ on $\cF$.
This is the point where mathematicians usually stop, since usually
such measures are ill-defined.
But let this not disturb us.

In the {\em quasi-classical limit} $k\to\infty$, the stationary
phase method (see e.g. \cite{stationary} and also Exercise
\ref{ex:stationary}) states that under some reasonable assumptions
about the behavior of $S$ this fast-oscillating integral localizes
on the critical points of $S$, so one recovers the classical case.

The {\em expectation value} $\<f\>$ of an {\em observable}
$f:\cF\to\R$ is
$$
   \<f\>=\frac1{Z}\int_{\cF}\cD\Gp\ e^{ikS(\Gp)}f(\Gp).
$$
For a collection $f_1$, \dots, $f_m$ of
observables their {\em correlation function} is
$$
   \<f_1,\dots,f_m\>= \frac1{Z}\int_{\cF}\cD\Gp\
       e^{ikS(\Gp)}\prod_{i=1}^n f_i(\Gp).
$$

By solving a theory one usually means a calculation of these
integrals or their asymptotics at $k\to\infty$.

Increasingly often, due to a simpler behavior and better convergence
properties, one considers instead the Euclidean partition function,
equally well encoding physical information (and related to
\eqref{eq:minkowski} by a certain analytic continuation in the time
domain, called Euclidean, or Wick, rotation):
\begin{equation}
\label{eq:euclidean}
Z=\int_{\cF}e^{-kS(\Gp)}\cD\Gp.
\end{equation}

Since at present
a general mathematical treatment of path integrals is lacking,
we will first consider a finite dimensional case.

\subsection{Finite-dimensional version of QFT}
\label{sub:fdim} Let us take $\cF=\R^d$ as the space of fields. An
action $S$ and observables $f_i$ are then just functions
$\R^d\to\R$. For a constant $k\in\R$, consider the partition
function $Z=\int_{\R^d}dx e^{-kS(x)}$ and the correlation
functions $\<f_1,\dots,f_m\>=Z^{-1}\int_{\R^d}dx\ e^{-kS(x)}\prod_i
f_i(x)$. We are interested in the behavior of $Z$ and
$\<f_1,\dots,f_m\>$ in the "quasi-classical limit" $k\to\infty$.

A well-known stationary phase method states that for large $k$ the
main contribution to $Z$ and $\<f_1,\dots,f_m\>$ comes from some
small neighborhoods of the points $x$ where $\dd S/\dd x=0$. Thus
it suffices to study a behavior of $Z$ and $\<f_1,\dots,f_m\>$
near such a point $x_0$. Considering the Taylor expansions of $S$
and $f_i$ in $x_0$ (and noticing that the linear terms in the
expansion of $S$ vanish), after an appropriate changes of
coordinates we arrive to the following problem: study integrals
$$\int_{\R^d}dx\ e^{-\frac12\<x,Ax\>+\h U(x)}P(x)$$ for some
bilinear form $A$, higher order terms $U(x)$, and monomials $P(x)$
in the coordinates $x^i$.

Further in these notes we will calculate such integrals
explicitly. To keep track of all terms appearing in these
calculations, we will use Feynman diagrams as a simple
book-keeping device. See the notes of Kazhdan in \cite{string}
for a more in-depth treatment.

\section{Finite-dimensional Feynman diagrams}
\label{sec:findim}

\subsection{Gauss integrals}
Recall a well-known formula for the Gauss integral (obtained by
calculating the square of this integral in polar coordinates):
\begin{prop}\label{prop:gauss}
$$
         \int_{-\infty}^\infty dx e^{-\frac12ax^2}=\sqrt{\frac{2\pi}a}.
$$
\end{prop}

More generally, let $A=(A_{ij})$ be a real $d\times d$
positive-definite matrix, $x=(x^1,\dots,x^d)$ the Euclidean
coordinates in $V=\R^d$, and $\<\ ,\ \>:(\R^d)^*\times\R^d\to\R$ the
standard pairing $\<x_i,x^j\>=\Gd_i^j$. Then

\begin{prop}
\begin{equation}\label{gauss}
Z_0=\int_{\R^d} dx\ e^{-\frac12\<Ax,x\>}=
\left(\det\frac{A}{2\pi}\right)^{-\frac12}.
\end{equation}
\end{prop}

Indeed, by an orthogonal transformation (which does not change the
integral) we can diagonalize $A$ and apply the previous formula in
each coordinate.

\begin{rem}
In a more formal setting, this may be considered as an equality
for a positive-definite symmetric operator $A:V\to V^*$ from a
$d$-dimensional vector space $V$ to its dual (and
$\<\ ,\ \>:V^*\times V\to\R$).
Indeed, $A$ induces $\det A:\GL^dV\to\GL^dV^*$, so that $\det
A\in(\GL^dV^*)^{\otimes2}$ and $(\det A)^{-\frac12}\in|\GL^dV|$.
Hence equality \eqref{gauss} with $\R^d$ changed to $V$ still
makes sense if we consider both sides as elements of $|\GL^dV|$.
In a similar way, for  $\C$-valued symmetric operator $A:V\to V^*$
with a positive-definite $\text{Im} A $ one has
$$  \displaystyle{\int_V dx\ e^{\frac{i}2\<Ax,x\>}=
    \left(\det\frac{A}{2\pi i}\right)^{-\frac12}}
$$
where now both sides belong to $|\GL^dV|_\C$.
\end{rem}

A more general form of equation \eqref{gauss} is obtained by
adding a linear term $-\<b,x\>$ with $b\in(\R^d)^*$ to the exponent:
define $Z_b$ by
\begin{equation}\label{def_Zb}
Z_b=\int dx\ e^{-\frac12\<Ax,x\>+\<b,x\>}.
\end{equation}
Then, by a change $x\to x-A^{-1}b$ of coordinates, we obtain
\begin{prop}
\begin{equation}\label{Z_b}
Z_b=\left(\det\frac{A}{2\pi}\right)^{-\frac12}e^{\frac12
\<b,A^{-1}b\>}=Z_0e^{\frac12\<b,A^{-1}b\>}.
\end{equation}
\end{prop}

\begin{exs}
\label{ex:stationary}
Verify the stationary phase method in the simplest case: use
an appropriate change of coordinates to pass
$$
               \text{from}\quad
   \int_{\Ga}^{\Gb} dx\ e^{k(-\frac{ax^2}{2}+bx)}
             \quad  \text{to}\quad
   \int_{\Ga'}^{\Gb'} dx\ e^{-\frac{x^2}{2}}.
$$
What happens to a small $\Ge$-neighborhood of the critical
point $x_0=b/a$ under this change of coordinates?
Conclude that in the limit $k\to\infty$ integration over a
small neighborhood of $x_0=b/a$ gives the same leading term
in the expansion of this integral in powers of $k$, as
integration over the whole of $\R$.
\end{exs}

\subsection{Correlation functions}
The correlators $\<f_1,\dots,f_m\>$ of $m$ functions
$f_i:\R^d\to\R$ (also called {\em $m$-point functions}) are
defined by plugging the product of these functions in the
integrand and normalizing:
\begin{equation}\label{def_corr}
\<f_1,f_2,\dots,f_m\>= \frac1{Z_0}\int dx\ e^{-\frac12\<Ax,x\>}
                        f_1(x)\dots f_m(x).
\end{equation}
They may be computed using $Z_b$. Indeed, notice that
$$
   \frac\dd{\dd b_i}  \int dx\ e^{-\frac12\<Ax,x\>+\<b,x\>}
        =\int dx\ e^{-\frac12\<Ax,x\>+\<b,x\>} x^i,
$$
hence for correlators
of any (not necessary distinct) coordinate functions we have
\begin{equation}\label{m_point}
  \<x^{i_1},\dots,x^{i_m}\>=\frac1{Z_0}\dd_{i_1} \dots\dd_{i_m}Z_b
\big|_{b=0}=\dd_{i_1} \dots \dd_{i_m}e^{\frac12\<b,A^{-1}b\>}
\big|_{b=0}
\end{equation}
where we denoted $\dd_i=\dd/ \dd b_i$.

In particular, 2-point functions are given by the Hessian matrix
$\frac{\dd^2}{\dd b^2}(Z_b/Z_0)\big|_{b=0}$ with the matrix
elements
\begin{equation}\label{eq:hess}
\<x^i,x^j\>=\dd_i\dd_j e^{\frac12 \<b,A^{-1}b\>}\big|_{b=0}
=(A^{-1})_{ij}.
\end{equation}
Thus the bilinear pairing $\Sym^2(V^*)\to\R$ given by 2-point
functions is just the pairing determined by $A^{-1}$. This
explains the similarity of our notations for the 2-point functions
and $\<\ ,\ \>:V^*\times V\to\R$.

For polynomials, or more generally, formal power series $f_1,\dots,
f_m$ in the coordinates we may apply \eqref{m_point} (with
$i_1=i_2=\dots=i_n=i$ for each monomial $(x^i)^n$) and then put
the series back together, noting that each $x^i$ should be
substituted by $\dd_i$. This yields:
\begin{prop}
\begin{equation}\label{corr}
\left. \<f_1,f_2,\dots,f_m\>=f_1 \left( \frac\dd{\dd b} \right)\dots
  f_m  \left( \frac\dd{\dd b}\right) \,  e^{\frac12\<b,A^{-1}b\>} \right|_{b=0}
\end{equation}
\end{prop}

\subsection{Wick's theorem}
Denote by $A^{ij}$ the matrix elements $(A^{-1})_{ij}$ of
$A^{-1}$. The key ingredient of the Feynman diagrams technique is
Wick's theorem (see e.g. \cite{Sh}) which we state in its simplest form:
\begin{thm}[Wick]
\begin{equation}\label{wick}
\dd_{i_1}\dots\dd_{i_m} e^{\frac12\<b,A^{-1}b\>}
\big|_{b=0}=\begin{cases}
\sum A^{j_1 j_2}\dots A^{j_{m-1} j_m},\quad m=2n\\
0,\quad m=2n+1\end{cases}
\end{equation}
where the sum is over all partitions $(j_1,j_2)$,\dots,
$(j_{m-1},j_m)$ in pairs of the set $i_1$,$i_2$,\dots,$i_m$ of
indices.
\end{thm}
\begin{proof}
For each $k$, the expression $\dd_{i_1}\dots\dd_{i_k} e^{\frac12\<b,A^{-1}b\>}$,
considered as a function of $b$, is always of the form
$P_{i_1\dots i_k}(b)e^{\frac12\<b,A^{-1}b\>}$, where
$P_{i_1\dots i_k}(b)$ is a polynomial.
Each new derivative $\dd_j$ acts either on the polynomial part,
or on the exponent, by the rule
$$
   \dd_j\big(P(b)e^{\frac12\<b,A^{-1}b\>}\big)=
    \dd_j(P(b))e^{\frac12\<b,A^{-1}b\>}+
    P(b)(\sum_i A^{ji}b_i)e^{\frac12\<b,A^{-1}b\>},
$$
so the polynomial part $P_{i_1\dots i_m}(b)$ may be defined
recursively by $P_{\emptyset}(b)={\bf 1}$ and
\begin{multline}
\label{eq:recursive}
P_{i_1\dots i_m}(b)=
(\dd_{i_1}+\sum_i A^{i_1 i}b_i)P_{i_2\dots i_m}(b)=\dots\\
=(\dd_{i_1}+\sum_i A^{i_1 i}b_i) \dots (\dd_{i_m}+\sum_i
A^{i_m i}b_i) {\bf 1},
\end{multline}
where ${\bf 1}$ is the function identically equal to 1.
We are interested in the constant term $P_{i_1\dots i_m}(0)$.
Directly from \eqref{eq:recursive} we can make two observations.
Firstly, if $m$ is odd, $P_{i_1\dots i_m}(b)$ contains only
terms of odd degrees, in particular $P_{i_1\dots i_m}(0)=0$.
Secondly, unless each derivative $\dd_{i_k}$, $k<m$ acts on the
term $\sum_iA^{i_l i}b_i$ in some $l$-th, $l>k$, factor of the
\eqref{eq:recursive}, the evaluation at $b=0$ would give zero.
Each such pair $(i_k,i_l)$ contributes a factor of $A^{i_li_k}$
to the constant  term of $P_{i_1\dots i_m}$.
These observations prove the theorem.
\end{proof}

It is convenient to extend \eqref{wick} by linearity to arbitrary
linear functions of the coordinates, note that in this case we may
define the $2$-point functions $\<f,g\>$ by $\<f,A^{-1}g\>$ in view
of \eqref{eq:hess}, and finally combine it with \eqref{m_point} into
the following version of Wick's theorem:

\begin{thm}[Wick]\label{thm:wick}
Let  $f_1(x),\dots,f_m(x)$ be arbitrary linear functions of the
coordinates $x_i$. Then all $m$-point functions vanish for odd $m$.
For $m=2n$ one has
\begin{equation}
\label{eq:symm}
\<f_1,\dots,f_m\>=\sum \<f_{i_1},f_{i_2}\>\dots\<f_{i_{m-1}},f_{i_m}\>,
\end{equation}
where the sum is over all pairings $(i_1,i_2)$,\dots, $(i_{m-1},i_m)$
of $1,\dots,m$ and the $2$-point functions $\<f_j,f_k\>$ are given by
$\<f_j,A^{-1}f_k\>$.

\end{thm}

\begin{rem}
Another idea for a proof of Theorem \ref{thm:wick} is the following.
Note that both sides of \eqref{eq:symm} are symmetric functions of
$1,\dots,m$, so they  may be considered as functions on $m$-th symmetric
power $S^m(V)$ of $V=\R^d$.
Thus it suffices to check \eqref{eq:symm} only for $f_1=\dots=f_m=f$;
in this case it is obvious.
\end{rem}

\begin{exs}
Check that the number of all pairings of $1,\dots,2n$ is
$(2n)! / 2^n n!$. Calculate
$$
   \int_{-\infty}^\infty dx\ x^m e^{x^2/2}
$$
 using integration by parts and Proposition~\ref{prop:gauss}. Calculate
$$
  \left. \frac{d^m}{d x^m} e^{x^2/2} \right|_{x=0}
$$
substituting $x^2/2$ instead of $x$ in the Taylor series expansion of $e^x$.
Compare these expressions and explain how are they related to the above number of pairings.
\end{exs}

\begin{exs}
Find formulas for the 4-point functions $\<x_1,x_1,x_2,x_3\>$ and
$\<x_1,x_1,x_1,x_2\>$.
\end{exs}

\subsection{First Feynman graphs}
It is convenient to represent each term
$$
   \<f_{i_1},f_{i_2}\>\dots\<f_{i_{m-1}},f_{i_m}\>
$$
in Wick's formula \eqref{eq:symm} by a simple graph. Indeed,
consider $m$ points, with the $k$-th point representing $f_k$.
A pairing of $1,\dots,2n$ gives a natural way to connect these
points by $n$ edges, with an edge (a {\em propagator} in the
physical jargon) $e=(j,k)$ representing
$A^{-1}_e=\<f_j,A^{-1}f_k\>$.
Equation \eqref{eq:symm} becomes then
\begin{equation}\label{m_diagr}
\<f_1,\dots,f_m\>=\sum_{\GG}\prod_{e\in\text{edges}
(\GG)}A^{-1}_e,
\end{equation}
where the sum is over all univalent graphs as above.

\begin{ex}\label{4_point}
An application of equation \eqref{m_diagr} for $n=2$ (see Figure
\ref{fig:pair}a) gives the following:
$$\begin{aligned}
\<x_1,x_2,x_3,x_4\>&=A^{12}A^{34}+A^{13}A^{24}+ A^{14}A^{23},\\
\<x_1,x_1,x_2,x_2\>&=A^{11}A^{22}+2A^{12}A^{12},\\
\<x_1,x_1,x_1,x_1\>&=3A^{11}A^{11}.
\end{aligned}$$
%
\end{ex}

\subsection{Adding a potential}
The above computations may be further generalized by adding a {\em
potential function} $U(x)$ (with some small parameter $\h=k^{-1}$)
to $\<Ax,x\>$ in the definition of $Z_0$.
Namely, define\footnote{Again, let me remind that we ignore problems
of convergence: for most $U(x)$ this integral will be divergent!}
$Z_U$ by
\begin{equation}\label{def_Z_U}
Z_U=\int dx\ e^{-\frac12\<Ax,x\>+\h U(x)}.
\end{equation}
Applying \eqref{corr} for $f=e^{\h U(x)}$ we get:

\begin{prop}
\begin{equation}\label{Z_U}
Z_U=Z_0 e^{\h U(\frac\dd{\dd
b})}e^{\frac12\<b,A^{-1}b\>}\big|_{b=0}.
\end{equation}
\end{prop}

Correlation functions $\<f_1,\dots,f_m\>_U$ are defined similarly
to \eqref{def_corr}:
\begin{equation}\label{def_corr_U}
\<f_1,f_2,\dots,f_m\>_U=\frac1{Z_U}\int dx\ e^{-\frac12\<Ax,x\>+\h
U(x)} f_1(x)\dots f_k(x).
\end{equation}
Using \eqref{corr} once again, we get

\begin{prop}
\begin{equation}\label{corr_U}
\<f_1,f_2,\dots,f_m\>_U=\frac{Z_0}{Z_U} e^{\h U(\frac\dd{\dd b})}
   f_1\left(\frac\dd{\dd b}\right)\dots f_m\left(\frac\dd{\dd b}\right)
           \left.  e^{\frac12\<b,A^{-1}b\>}\right|_{b=0}.
\end{equation}
\end{prop}

\subsection{A cubic potential}\label{sub:cubic}
\begin{figure}
\includegraphics[height=0.6in]{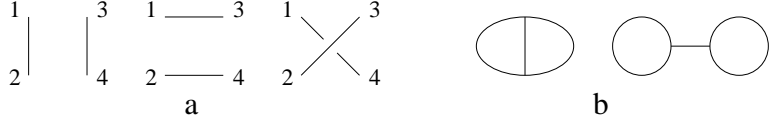}
\caption{Terms of $\<x_1,x_2,x_3,x_4\>$ and graphs of degree two}
\label{fig:pair}
\end{figure}

Consider the important example of a cubic potential function
                 $U(x)= \sum U_{ijk}x^i x^j x^k$.
Let us compute the expansion of the partition function \eqref{def_Z_U}
in power series in $\h$.
The coefficient of $\h^n$ in the expansion of \eqref{Z_U} is
$$
  \frac{Z_0}{n!} \left( \sum_{i,j,k} U_{ijk}\dd_i\dd_j\dd_k \right)^n\,
     \left.  e^{\frac12\<b,A^{-1}b\>}\right|_{b=0}.
$$
Let us start with the lowest
degrees. By Wick's theorem, the coefficient of $\h$ vanishes and
the coefficient of $\h^2$ is given by
\begin{multline}\label{h2}
\frac{Z_0}{2!}\sum_{i,j,k}\sum_{i',j',k'} U_{ijk}U_{i'j'k'}\dd_i
\dd_j\dd_k\dd_{i'}\dd_{j'}\dd_{k'}
e^{\frac12\<b,A^{-1}b\>} \big|_{b=0} =\\
\frac{Z_0}{2!}\sum_{i,j,k}\sum_{i',j',k'} U_{i j k} U_{i'j'k'}
              \sum  A^{i_1i_2}A^{i_3i_4}A^{i_5i_6},
\end{multline}
where the last sum is over all pairings
$(i_1,i_2)$,\dots,$(i_5,i_6)$ of $i,j,k,i',j',k'$. We may again
encode these pairings by labelled graphs, connecting 6 vertices
labelled by $i,j,k,i',j',k'$ by three edges
$(i_1,i_2)$,$(i_3,i_4)$,$(i_5,i_6)$ representing
$A^{i_1i_2}A^{i_3i_4}A^{i_5i_6}$. This time, however, we have an
additional factor $U_{i j k}U_{i'j'k'}$. To represent $U_{i j k}$
graphically, let us glue the triple $(i,j,k)$ of univalent
vertices in a trivalent vertex; to preserve the labels, we can
write them on the ends of the edges meeting in this new vertex
(i.e., on the star of the vertex). Similarly, we represent
$U_{i'j'k'}$ by gluing the remaining triple $(i',j',k')$ of
univalent vertices into a second trivalent vertex. Thus for each
of the $6!/ (2^3 3!) =15$ pairings of $i,j,k,i',j',k'$ we end
up with a graph with two trivalent vertices; we get 6 copies of
the $\GT$-graph and 9 copies of the dumbbell graph shown in Figure
\ref{fig:pair}b.

Note, however, that each of these labelled graphs is considered
up to its automorphisms, i.e. maps of a graph onto itself, mapping
edges to edges and vertices to vertices and preserving the incidence
relation.
Indeed, while the application of an automorphism changes the labels,
it preserves their pairing (edges) and the way they are united in
triples (vertices), thus corresponds to the same term in the right
hand side of \eqref{h2}.
Instead of summing over the automorphism classes of graphs, we may
sum over all labelled graphs, but divide the term corresponding to
a graph $\GG$ by the number $|\Aut \GG |$ of its automorphisms.
E.g., for the $\GT$-graph of Figure \ref{fig:pair}b $|\Aut \GG |=12$,
and twelve copies of this graph (which differ only by transpositions
of the labels) all give the same terms
$U_{ijk}U_{i'j'k'}A^{ii'}A^{jj'}A^{kk'}$.

Also, when summing the resulting expressions over all indices,
note that the terms corresponding to $i,j,k,i',j',k'$ and to
$i',j',k',i,j,k$ are the same (which will cancel out with $1/2!$
in front of the sum). Hence, we may write the coefficient of
$\h^2$ in the following form:
$$Z_0\sum_\GG\frac1{|\Aut \GG|}\sum_{labels}
\prod_{v} U_v\prod_{e} A^{-1}_e,$$ Here the sum is over all
trivalent graphs with two vertices and labellings of their edges,
$U_v=U_{ijk}$ for a vertex $v$ with the labels $i,j,k$ of the
adjacent edges, and $A^{-1}_e=A^{ij}$ for an edge $e$ with labels
$i,j$.

\begin{exs}
Calculate the number of automorphisms of the dumbbell graph of
Figure \ref{fig:pair}b.
\end{exs}

In general, for the coefficient of $\h^n$ we get the same formula,
but with the summation being over all labelled trivalent graphs
with $n$ vertices.

\subsection{Correlators for a cubic potential}
\begin{figure}
\includegraphics[width=5in]{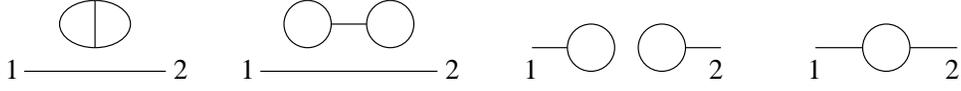}
\caption{Degree two graphs with two legs} \label{fig_corr}
\end{figure}

We may treat $m$-point functions in a similar way. Let us first
consider the power series expansion in $\h$ of
$Z_U\<x^{i_1},\dots,x^{i_m}\>_U$. The coefficient of $\h^n$ is
$$
     \frac{Z_0}{n!}\left(\sum_{i,j,k} U_{ijk}\dd_i\dd_j\dd_k\right)^n
     \dd_{i_1}\dots\dd_{i_m}\,
       \left.  e^{\frac12\<b,A^{-1}b\>}\right|_{b=0}.
$$
Thus it may again be presented by a sum over labelled graphs, with
the only difference being that now in addition to $n$ trivalent
vertices these graphs also have $m$ ordered {\em legs} (i.e.
univalent vertices) labelled by $i_1,\dots,i_m$. See Figure
\ref{fig_corr} for graphs representing the coefficient of $\h^2$
in $Z_U\<x^1,x^2\>_U$.

However, not all of these graphs will enter in the expression for
$\<x^{i_1},\dots,x^{i_m}\>_U$, since we should now divide this sum
over graphs by $Z_U$ (represented by a similar sum, but over
graphs with no legs). This will remove all {\em vacuum diagrams},
i.e. all graphs which contain some component with no legs. Indeed,
the term corresponding to a non-connected graph is a product of
terms corresponding to each connected component. Each component
with no legs appears also in the expansion of $Z_U$ and thus will
cancel out after we divide by $Z_U$. For example, the first graph
of Figure \ref{fig_corr} contains a vacuum $\GT$-graph component.
But it also appears in the expansion of $Z_U$ (see Figure
\ref{fig:pair}b). Thus the corresponding factor cancels out after
division by $Z_U$. The same happens with the second graph of
Figure \ref{fig_corr}. As a result, only the last two graphs of
Figure \ref{fig_corr} will contribute to the coefficient of $\h^2$
in the expansion of $\<x^1, x^2\>_U$.

\begin{ex}["A finite dimensional $\phi^3$-theory"]
Take $U_{ijk}=\Gd_{ij}\Gd_{jk}$, i.e. $U=\sum_i(x^i)^3$. Note that
since all ends of edges meeting in a vertex are labelled by the
same index, we may instead label the vertices. Thus the
calculation rules are quite simple: we count uni-trivalent graphs;
a vertex represents a sum $\sum_i$ over its labels; an edge with
the ends labelled by $i,j$ represents $A^{ij}$. The coefficients
of $\h^2$ in $Z_U$ and in $\<x^1,x^2\>_U$ are given by
$$
    \begin{aligned}
Z_0&\sum_{i,j} 6(A^{ij})^3+9A^{ij}A^{ii}A^{jj},\\&
\sum_{i,j}9A^{1i}A^{2j}A^{ii}A^{jj}+6A^{1i}A^{2j}(A^{ij})^2
     \end{aligned}
$$
respectively. We can identify these terms with two graphs of
Figure \ref{fig:pair}b and two last graphs of Figure
\ref{fig_corr}, respectively. The first two graphs of Figure
\ref{fig_corr} represent
$$
   \sum_{i,j} 6A^{12}(A^{ij})^3+9A^{12}A^{ij}A^{ii}A^{jj}.
$$
These terms do appear in the $\h^2$ coefficient of
$Z_U\<x^1,x^2\>$, but cancel out after we divide it by
$Z_U=Z_0\big(1+\h^2\sum_{i,j}(6(A^{ij})^3+9A^{ij}A^{ii}A^{jj})
+\dots\big)$.
\end{ex}

\subsection{General Feynman graphs}
It is now clear how to generalize the above results to the case of
a general potential $U(x)$: the $k$-th degree term $U_{i_1\dots
i_k}x^{i_1}\dots x^{i_k}$ of $U$ will lead to an appearance of
$k$-valent vertices representing factors $U_{i_1\dots i_k}$. We
will call such a vertex an {\em internal vertex}. We assume that
there are no linear and quadratic terms in the potential, so
further we will always assume that all internal vertices of any
Feynman graph $\GG$ are of valence $\ge 3$; denote their number by
$|\GG|$. Denote by $\GG^0$ the set of all graphs with no legs.
Also, for $m\ge1$, denote by $\GG^m$ the set of all non-vacuum
(i.e. such that each connected component has at least one leg)
graphs with $m$ ordered legs.

Denoting $U_v=U_{i_1\dots i_k}$ for an internal vertex $v$ with
the labels $i_1,\dots,i_k$ of adjacent edges, and
$A^{-1}_e=A^{ij}$ for an edge $e$ with its ends labelled by $i,j$,
we get
\begin{prop}
\begin{equation}
\label{eq:ZU} Z_U=Z_0\sum_{\GG\in\GG^0}\frac{\h^{|\GG|}}{|\Aut
\GG|}\sum_{labels} \prod_{v} U_v \prod_{e} A^{-1}_e .
\end{equation}
\end{prop}
Note that instead of performing the internal summation over all
labellings, one may include the summation over labels of the star
of a vertex into the weight of this vertex.

In a similar way, for $m$-point functions we get
\begin{prop} For even $m$,
\begin{equation}
\label{eq:corr}
\<x^{i_1},\dots,x^{i_m}\>_U=\sum_{\GG\in\GG^m}\frac{\h^{|\GG|}}
{|\Aut \GG|}\sum_{labels}\prod_{v} U_v\prod_{e}A^{-1}_e ,
\end{equation}
where the sum is over all labelled graphs $\GG$ with $m$ legs
labelled by $i_1,\dots,i_m$.
\end{prop}
Again, we may include the summation over the labels of the star of
an internal vertex into the weight of this vertex.

\subsection{Weights of graphs}
Let us reformulate the above results using a general notion of
weights of graphs.

Let $V$ be a vector space. A {\em weight system} is a collection
$(a,\{u_k\}_{k=3}^\infty)$ of $a\in\Sym^2(V)$ and
$u_k\in\Sym^k(V^*)$. A weight system $W$ defines a {\em weight}
$W_\GG:(V^*)^{\otimes m}\to\R$ of a graph $\GG\in\GG^m$ in the
following way. Assign $u_k\in\Sym^k(V^*)$ to each internal vertex
$v$ of valence $k$, associating each copy of $V^*$ with (an end
of) an edge. Also, to the $i$-th leg of $\GG$, $i=1,\dots,m$
assign some $f_i\in V^*$. Now, for each edge contract two copies
of $V^*$ associated to its ends using $a\in\Sym^2(V)$. After all
copies of $V^*$ get contracted, we obtain a number
$W_\GG(f_1,\dots,f_m)\in\R$.

In our case, a bilinear form $A^{-1}$ and a potential $\h U(x)$
determine a weight system in an obvious way: set $a=A^{-1}$ and
let $u_v$ to be the degree $k$ part of $\h U(x)$.
These rules of computing the weights corresponding to a physical theory are
called {\em Feynman rules}.

Formulas \eqref{eq:ZU} and \eqref{eq:corr} above can be
reformulated in these terms as
\begin{equation}
\label{eq:graphs}
\begin{aligned}
             Z_U &=Z_0\sum_{\GG\in\GG^0}\frac1{|\Aut \GG|}W_\GG, \\
    \<f_1,\dots,f_m\>_U &=\sum_{\GG\in \GG^m} \frac1{|\Aut \GG|}
                      W_\GG(f_1,\dots,f_m).
\end{aligned}
\end{equation}

\begin{exs}[Finite dimensional $\phi^4$-theory]
Consider a potential $U=\sum_i (x^i)^4$. Formulate the Feynman
rules. Find the graphs which contribute to the coefficient of
$\h^2$ of $Z_U$ and compute their coefficients. Do the same for
$\<x^1,x^2\>_U$. Draw the graph representing $\sum_i
A^{12}(A^{ii})^2$; does it appear in the expansion of
$\<x^1,x^2\>_U$ and why?
\end{exs}

\subsection{Free energy: taking the logarithm}
The summation in equation \eqref{eq:graphs} is over all graphs in
$\GG^0$, which are plenty. Denote by $\GG^0_{conn}$ the subset of
all connected graphs in $\GG^0$. There is a simple way to leave
only a sum over graphs in $\GG^0_{conn}$, namely to take the
logarithm of the partition function (called the {\em free energy}
in the physical literature):

\begin{prop}\label{log}
Let $W$ be a weight system. Then
$$
    \log\left (\sum_{\GG\in\GG^0}\frac1{|\Aut\GG|}W_\GG \right)=
      \sum_{\GG\in\GG^0_{conn}}\frac1{|\Aut \GG|}W_\GG .
$$
\end{prop}

\begin{proof}
Let us compare the terms of the power series expansion for the
right hand side with the terms in the left hand side:
$$
  \exp\left(\sum_{\GG\in\GG^0_{conn}}\frac1{|\Aut \GG|}W_\GG\right)=
   \sum\frac1{n_1!\dots n_k!} W_{\GG_1}^{n_1}\dots W_{\GG_k}^{n_k},
$$
where the sum is over all $k$, $n_i$, and distinct
$\GG_i\in\GG^0_{conn}$, $i=1,\dots,k$.
Consider $\GG=(\GG_1)^{n_1}\dots(\GG_k)^{n_k}\in\GG^0$.
Since in addition to automorphisms of each $\GG_i$ there are also
automorphisms of $\GG$ interchanging the $n_i$ copies of $\GG_i$,
we have $|\Aut\GG|=n_1!\dots n_k!|\Aut\GG_1|\dots|\Aut\GG_k|$.
Also, any weight system satisfies $W_{\GG'\GG''}=W_{\GG'}W_{\GG''}$,
hence $W_\GG= W_{\GG_1}^{n_1}\dots W_{\GG_k}^{n_k}$.
The proposition follows.
\end{proof}

\begin{exs}
Formulate and prove a similar statement for graphs with legs.
\end{exs}

\begin{rem}
It is possible to restrict the class of graphs to 1-connected (in
the physical literature usually called {\em $1$-point
irreducible,} or 1PI for short) graphs. A graph is 1-connected, if
it remains connected after a removal of any one of its edges. This
involves a passage to a so-called effective action, which I will
not discuss here in details. Mathematically, it simply means an
application of a Legendrian transform (a discrete version of a
Fourier transform): if $z(b)=\log(Z_b)$ is given by the sum over
all connected graphs as in Proposition \ref{log}, then
$\hat{z}(x)=\<b,x\>-z(b)$ is given by a similar sum over all 1PI
graphs (and $b(x)$ may be recovered as $\dd\hat{z}/\dd x$).
\end{rem}
\section{Gauge theories and gauge fixing}
\subsection{Gauge fixing}
All calculations of the previous section dealt only with the case
of a non-degenerate bilinear form $A$; in particular, the critical
points of the action $S(x)$ had to be isolated (see Section
\ref{sub:fdim}). However, gauge theories present a large class of
examples when it is not so. Suppose that we have an
$l$-dimensional group of symmetries, i.e. the Lagrangian is
invariant under a (free, proper, isometric) action of an $l$-dimensional
Lie group $G$. Then instead of isolated critical points we have
critical orbits, so $A$ has $l$ degenerate directions and the
technique of Gauss integration can not be applied.

Let us try to calculate the partition and correlation functions
without a superfluous integration over the orbits of $G$.
In other words, we wish to reduce integrals of $G$-invariant
functions on $X$ to integrals on the quotient space $\tX=X/G$
of $G$-orbits.
For this purpose, starting from a $G$-invariant measure on $X$
we should desintegrate it as the Haar measure  on the orbits over
some ``quotient measure''  $\tilde\mu$ on the base $\tX$.

If $G$ is compact then $\tilde\mu$ is the standard push-forward of $\mu$.
For example, if $f$ is a rotationally invariant function on $\R^2$,
we can take the pair of polar coordinates $(r,\phi)$ as coordinates
in the quotient space $\tX$ and the orbit, respectively.
The measure on $\tX$ in this case is $2\pi r\, dr$ and we get the
following elementary formula:
$$
  \int_{\R^2}f(|x|)\ d^2x=2\pi\int_0^\infty f(r) r\ dr.
$$

If $G$ is a locally compact group acting properly on $X$ then $\tilde\mu$ can be
defined by the property
$$
    \mu(Y) = \int_{\tX} | Y\cap G(\tilde x)|\, d\tilde\mu(\tilde x),
$$
where $G(\tilde x)$ is the fiber over $\tilde x\in \tX$ and $|\cdot|$ is the Haar measure on it.
In this case the integral $\int_X f\, dx$  in question is infinite,
but it can be formally defined  (``regularized'') as $\int_{\tilde X} f d\tilde x$.

A standard physical procedure for the desintegration that can be
applied also to a non-locally compact gauge group is called a
{\em gauge  fixing} (see e.g., \cite{Sc}); it goes as follows.
Suppose that $f:X\to\R$ is $G$-invariant, i.e. $f(gx)=f(x)$ for
all $x\in X$, $g\in G$. Choose a (local) section $s:\tX\to
X$ which intersects each orbit of $G$ exactly once.
Suppose that it is defined by $l$ independent equations
$F^1(x)=\dots=F^l(x)=0$ for some $F:X\to\R^l$.
Firstly, we want to count each $G$-orbit only once.
This is simple to arrange by inserting an $l$-dimensional
$\Gd$-function $\Gd^l(F(x))$ in the integrand. Secondly, we want
to take into account the volume of a $G$-orbit passing through $x$,
so we should count each orbit with a certain Jacobian factor
$J(x)$ (called the {\em Faddeev-Popov determinant}).
How should one define such a factor? We wish to have
$$
  \int_X f(x)\ dx=\int_X f(x) J(x) \Gd^l(F(x))\ dx.
$$
Rewriting the right hand side to include an additional
integration over $G$ and noticing that both $f(x)$ and
$J(x)$ are $G$-invariant, we get
\begin{multline*}
    \int_X f(x)\ dx=\int_X f(x) J(x) \Gd^l(F(x)\ dx=\\
    =\int_X dx \int_G dg\ f(x) J(x) \Gd^l(F(gx))=
    \int_X dx\ f(x) J(x)\int_Gdg\ \Gd^l(F(gx)).
\end{multline*}
Thus we see that we should define $J(x)$ by
$$
           J(x)\int_G dg\ \Gd^l(F(gx))=1,
$$
where $dg$ is the left $G$-invariant measure on $G$.
Thus, the Faddeev-Popov determinant plays the
role of Jacobian for a change of coordinates from $x$ to
$(s(\tilde{x}),g)$.
Example in \S \ref{ex:CP1} below provides a good illustration.

\begin{rem}
A formal coordinate-free way to define $J(x)$ is as follows.
The section $s:\tX\to X$ determines a push-forward
$s_*:T_{\tilde{x}}\tX \to\ T_x X$ of the tangent spaces.
The tangent space $T_x X$ thus decomposes as
$s_*\oplus i:T_{\tilde{x}}\tX\oplus\fG \to T_x X$, where
$i:\fG\to T_x X$ is the tangent space to the orbit,
generated by the Lie algebra $\fG$ of $G$.
The Jacobian $J(x)$ may be then defined as
$J(x)=\det(s_*\oplus i)$.
\end{rem}

\begin{rem}
Equivalently, one may note that the tangent space to the fiber at
$x\in s$ may be identified with $\fG$, to directly set
$J(x)=\det \GL $, where $\GL=(\frac{\dd F^i}{\dd \fG^j})$ and
$\{\fG^j\}_{j=1}^l$ is a set of generators of the Lie algebra
$\fG$ of $G$, see e.g. \cite{BN}.
I.e., $J(x)$ is the inverse ratio of the volume element of
$\fG$ and its image in $\R^l$ under the action of $G$ composed
with $F$.

Indeed, since $F$ has a unique zero on each orbit and since (due
to the presence of the delta-function) we integrate only near the
section $s$, we can use $F$ as a local coordinate in the fiber
over $x$. Making a formal change of variables from $g$ to $F$ we get
$$
J(x)^{-1}=\int_Gdg\ \Gd^l(F(gx))=\int_GdF\ \Gd^l(F(gx))\
\left. \det\left(\frac{\dd g}{\dd F}\right)=\det\left(\frac{\dd g}{\dd F}\right)\right|_{F=0}.
$$
Calculating $(\frac{\dd F(gx)}{\dd g})\big|_{F=0}$ at a point $x\in s$
and identifying the tangent space to the fiber with $\fG$, we obtain
$(\frac{\dd F}{\dd g})\big|_{F=0}=(\frac{\dd F^i}{\dd \fG^j})$.
\end{rem}

\begin{exs}
Let us return to the simple example of a rotationally invariant
function $f(x_1,x_2)=f(|x|)$ on $\R^2$, using this time the
gauge-fixing procedure.
The group $G=S^1$ acts by rotations: $\phi x=e^{i\phi}x$ and the
(normalized) measure on $G$ is $\frac1{2\pi}d\phi$.
We should use the positive $x_1$-axis for a section $s$, so we
may take e.g. $F=x_2$.
A slight complication is that the equation $x_2=0$ defines the
whole $x_1$-axis and not only its positive half, so each fiber of
$G$ intersects it twice and not once. This can be taken care of,
either by dividing the resulting gauge-fixed integral by two, or
by restricting its domain of integration to the right half-plane
$\R_+^2$ in $\R^2$.
In any case, using $x_2$ instead of $\phi$ as a local coordinate
in the fiber $Gx$ near $x\in s$ we get
$d\phi=d\big(\arctan(x_2/x_1)\big)=x_1 |x|^{-2} dx_2$
Thus for $x\in s$ we have
$$
J(x)^{-1}=\int \Gd(F(\phi x))\, \frac1{2\pi}d\phi=
\frac1{2\pi}\int \Gd(x_2) x_1 |x|^{-2}dx_2=\frac1{2\pi x_1} ,
$$
so $J(|x|)=2\pi |x|$ as expected and
$$
\int_{\R^2}f(|x|)\, d^2x=\int_{\R_+^2}f(|x|)2\pi|x|\Gd(x_2)\,
dx_1 dx_2=2\pi\int_0^\infty f(r) r\, dr .
$$
\end{exs}

\subsection{Faddeev-Popov ghosts}
After performing the gauge-fixing, we are left with the gauged-fixed
partition function
$$
        Z_{GF}=\int_{\R^d}dx\ e^{-\frac12\<Ax,x\>}\Gd^l(F(x))\det \GL .
$$
We would like to make it into an integral of the type we have been
studying before. We have two problems: to include
$\Gd(F(x))\det \GL$ in the exponent (i.e., in the Lagrangian)
and--- more importantly--- to make $A$ into a non-degenerate bilinear
form.

The $\Gd$-function is easy to write as an exponent using the
Fourier transform:
$$
     \Gd^l(F(x))=(2\pi)^{-l}\int_{\R^l}d\xi\ e^{i\<\xi,F(x)\>}.
$$
The gauge variables $\xi$ (called {\em Lagrange multipliers})
supplement the variables $x$, and the quadratic part of $\<\xi,F(x)\>$
supplements $\<Ax,x\>$ so that the quadratic part $A_F$ of the
gauge-fixed Lagrangian is non-degenerate.

The $\det\GL$ term is somewhat more complicated; it can be also
represented as a Gaussian integral, but over {\em anti-commuting}
variables $c=(c^1,\dots,c^l)$ and  $\bar{c}=(\bar{c}^1,\dots,\bar{c}^l)$,
called  {\em Faddeev-Popov ghosts}.
Thus
\begin{equation}\label{eq:commut}
c^ic^j+c^jc^i=\bar{c}^ic^j+c^j\bar{c}^i=\bar{c}^i\bar{c}^j
+\bar{c}^j\bar{c}^i=0
\end{equation}
There are standard rules of integration over anti-commuting
variables (known to mathematicians as the {\em Berezin integral}\/,
see e.g. \cite[Chapter 33]{Sc} and \cite{FS,IZ}).
The ones relevant for us are

$$\mbox{
$\int c^i\ dc^j=\int \bar{c}^i\ d\bar{c}^j=\Gd^{ij}$ and
$\int 1\ dc^j=\int 1\ d\bar{c}^j=0$.
}$$
The multiple integration (over e.g., $dc=dc^l\dots dc^1$)  is
defined by iteration. One may show that this implies (see the
Exercise below) that for any matrix $\GL$
$$
  \int e^{\<\bar{c},\GL c\>}\ dc\, d\bar{c}=\det \GL .
$$

\begin{exs}
Let $l=1$ and define the exponent $e^{\Gl\bar{c}c}$ by the
corresponding power series. Use the commutation relations
\eqref{eq:commut} to verify that only the two first terms of
this expansion do not vanish. Now, use the integration rules
to deduce that $\int e^{\Gl\bar{c}c}\ dc\, d\bar{c}=
\int(1+\Gl\bar{c}c)\ dc\, d\bar{c}=\Gl$.
\end{exs}

Thus we may rewrite $Z_{GF}$ by adding to the Lagrangian the
gauge-fixing term and the ghost term:
$$
Z_{GF}=\int dx\, d\xi\, dc\, d\bar{c}\ e^{-\frac12\<Ax,x\>+
\<\bar{c},\GL c\>+i\<\xi,F(x)\>} .
$$

At this stage we may again apply the Feynman diagram expansion to
the gauge-fixed Lagrangian. The Feynman rules change in an obvious
fashion. The quadratic form now consists of two parts: $A_F$ and
$\GL$, so there are two types of edges. The first type presents
$A_F$, with the labels $x^i$ and $\xi^i$ at the ends. The second
type presents $\GL$, with the labels $c^i$ and $\bar{c}^i$ at the
ends. Note that since $\GL$ is not symmetric, these edges are {\em
directed}. Also, there are new vertices, presenting all higher
degree terms of the Lagrangian (in particular some where edges of
both types meet).
An example of the Chern-Simons theory will be provided in Section
\ref{sec:CS}.

\subsection{An example of gauge-fixing}\label{ex:CP1}
Let us illustrate the idea of gauge-fixing on an example of the
standard $\C^*$-action on $\C^2$. In the coordinates
$(x_1,\bx_1,x_2,\bx_2)$ on $\C^2$ the gauge group acts by
$x_i\to\Gl x_i$, $\bx_i\to\bar{\Gl}\bx_i$. Let us take
$A=\frac{x_1}{x_2}\frac{\bx_1}{\bx_2}$ as an invariant
function.

Of course, the orbit space $\C P^1$ is quite simple and an
appropriate measure on $\C P^1$ is well known; in the coordinates
$z=x_1/x_2$, $\bar{z}=\bx_1/\bx_2$ it is given by $dz
d\bar{z}/(1+z\bar{z})^2$ We are thus interested in
\begin{equation}
\label{eq:direct} Z_{GF}=\int\frac{dz d\bar{z}}{(1+z\bar{z})^2}
e^{-\frac12 z\bar{z}} .
\end{equation}
Let us pretend, however, that we do not know this and proceed with
the gauge-fixing method instead.

The invariant measure on $\C^2$ is $dx_1 d\bx_1 dx_2
d\bx_2/(x_1\bx_1+x_2\bx_2)^2$. In a gauge $F=0$ we have
$$
 Z_{GF}=\int\frac{dx_1 dx_2 d\bx_1 d\bx_2}{(x_1\bx_1+x_2\bx_2)^2}
 e^{-\frac12\frac{x_1}{x_2}\frac{\bx_1}{\bx_2}}\Gd^2(F(x,\bx))
\det \GL ,$$ $$\text{where}\quad \GL=\left|\begin{array}{cc}
x_1F_{x_1}+x_2F_{x_2}&\bx_1F_{\bx_1}+\bx_2F_{\bx_2}\\
x_1\bar{F}_{x_1}+x_2\bar{F}_{x_2}&\bx_1\bar{F}_{\bx_1}+
\bx_2\bar{F}_{\bx_2}\end{array}\right|.$$ E.g., for $F=x_2-1$ we
get $\Gd^2(|x_2-1|)$ and $\det \GL =x_2\bx_2$.

\begin{exs}[Different gauges give the same result]
Consider $F=x_2^\Ga-1$. Show that $\Gd^2(|x_2^\Ga-1|)=|\Ga
x_2^{\Ga-1}|^{-2}\Gd^2(|x_2-1|)$ and
$\det \GL =\Ga\bar{\Ga}(x_2\bx_2)^\Ga$. Check that the dependence
on $\Ga$ in $Z_{GF}$ cancels out, thus gives the same result as
$F=x_2-1$. Show that it coincides with formula \eqref{eq:direct}.
\end{exs}

Finally, let us check that while the initial quadratic form $A$ is
degenerate, the supplemented quadratic form $A_F$ is indeed
non-degenerate. It is convenient to make a coordinate change
$x_1'=x_1$, $x_2'=x_2-1$. Using a Fourier transform we get
$$
   \Gd^2(x_2-1)=(2\pi)^{-2}\int d\xi d\bar{\xi}
        e^{i(\xi x_2'-\bar{\xi}\bx_2')}.
$$
Also, we have
$\frac{x_1}{x_2}=x_1'+x_1'\sum_{n=1}^\infty (-1)^n x_2'$. We can
now compute $A$ and $A_F$; in the coordinates
$(x_1',\bx_1',x_2',\bx_2')$ and
$(x_1',\bx_1',x_2',\bx_2',\xi,\bar{\xi})$, respectively, we have:
$$
  A=\left|\begin{array}{cccc}0&1&0&0\\1&0&0&0\\0&0&0&0\\0&0&0&0
   \end{array}\right| ,
\quad A_F=\left|\begin{array}{cccccc}0&1&0&0&0&0\\
1&0&0&0&0&0\\0&0&0&0&i&0\\0&0&0&0&0&-i\\0&0&i&0&0&0\\0&0&0&-i&0&0
\end{array}\right| .
$$

\section{Infinite dimensional case}

\subsection{The dictionary}
Path integrals are generally badly defined, so instead of trying
to deduce the relevant results rigorously, we will just provide a
basic dictionary to translate the finite dimensional results to
the infinite dimensional case.

The main change is that instead of the discrete set
$i\in\{1,\dots,d\}$ of indices we now have a continuous variable
$x\in M^n$ (say, in $\R^n$), so we have to change all related
notions accordingly. The sum over $i$ becomes an integral over
$x$. Vectors $x=x(i)=(x^1,\dots,x^d)$ and $b=b(i)$ become fields
$\phi=\phi(x)$ and $J(x)$. A quadratic form $A=A(i,j)$ becomes an
integral kernel $K=K(x,y)$.
Pairings $\<Ax,x\>=\sum_{i,j}x^i A_{ij} x^j$ and $\<b,x\>=\sum_i
b^ix^i$ become $\<K\phi,\phi\>=\int dxdy\ \phi(x)K(x,y)\phi(y)$
and $\<J,\phi\>=\int dx\ J(x)\phi(x)$ respectively.
The partition function $Z_b$ defined by \eqref{def_Zb} becomes a
path integral $Z_J$ over the space $\cF$ of fields
$$
   Z_J=\int \cD\phi\ e^{-\frac12\<K\phi,\phi\>+\<J,\phi\>}.
$$

The inverse $A^{-1}$ of $A$ defined by
   $\sum_k A_{ik}A^{kj}= \Gd_i^j$
corresponds now to the inverse $G=K^{-1}$ of $K$ defined by
$$
        \int dz\ K(x,z)G(z,y)=\Gd(x-y).
$$
Formula \eqref{Z_b} for $Z_b$ then translates into
$$
         Z_J=Z_0 e^{\frac12\<J,G J\>}.
$$
Correlators $\<x^{i_1},\dots,x^{i_m}\>$ defined by
\eqref{def_corr} become now $m$-point functions
$$
     \<\phi(x_1),\dots,\phi(x_m)\>=\frac1{Z_0}\int \cD\phi\
     e^{-1/2\<K\phi,\phi\>}\phi(x_1)\dots\phi(x_m).
$$

\subsection{Functional derivation}
A counterpart of the derivatives $\dd/ \dd x^i$ is given by
the {\em functional derivatives} $\Gd / \Gd \phi(x) $.
The theory of functional derivation is well-presented in many places
(see e.g. \cite{DFN}), so I will just briefly recall the main notions.
Let $F(\phi)$ be a functional. If the differential
$$
     DF(\phi) (\rho) = \lim_{\Ge\to 0}\frac{F(\phi+\Ge\rho)-F(\phi)}{\Ge}
$$
can be represented as
$\int \rho(x) h(x) dx$ for some function $h(x)$, then we
define $\frac{\Gd F}{\Gd\phi(x)}=h(x)$.
In general, the functional derivative $\Gd F/ \Gd\phi(x)$ is the
distribution representing the differential of $F$ at
$\phi$.
The reader can entertain himself by making sense of the following
formulas, which show that its properties are similar to usual
derivatives:
$$
         \frac{\Gd}{\Gd\phi(x)}\phi(y)=\Gd(x-y),
$$
$$
   \frac{\Gd}{\Gd\phi(x)}(F(\phi)H(\phi))=\frac{\Gd}{\Gd\phi(x)}(F(\phi))
   \cdot H(\phi)+F(\phi)\cdot\frac{\Gd}{\Gd\phi(x)}(H(\phi)).
$$

\begin{ex}
\begin{multline*}
\frac{\Gd}{\Gd J(y)}e^{\<J,\phi\>}=\frac{\Gd}{\Gd J(y)} e^{\int dx
J(x)\phi(x)}=\phi(y)e^{\int dx J(x)\phi(x)} =\phi(y)
e^{\<J,\phi\>} .
\end{multline*}
\end{ex}
\begin{exs}
Consider a (symmetric) potential function
\begin{equation}
\label{U} U(\phi)=\sum_n\frac1{n!}\int dx_1\dots dx_n\
U_n(x_1,\dots,x_n) \phi(x_1)\dots\phi(x_n).
\end{equation}
Prove that
$$
   \frac{\Gd}{\Gd\phi(y)}U(\phi)=\sum_n\frac1{n!}\int dx_1\dots dx_n\
         U_{n+1}(y,x_1,\dots,x_n)\phi(x_1)\dots\phi(x_n).
$$
\end{exs}
The inverse $G(x,y)$ can be written as a Hessian, similarly to
equation \eqref{eq:hess} for $A^{-1}$:
$$
   G(x,y)=\left. \frac1{Z_0}\frac{\Gd}{\Gd J(x)}\frac{\Gd}{\Gd J(y)}Z_J\right|_{J=0}.
$$
 More generally, for $m$-point functions we have, similarly to \eqref{corr},
$$
  \left.   \<\phi(x_1),\dots,\phi(x_m)\>=\frac1{Z_0}\frac{\Gd}{\Gd J(x_1)}
            \dots\frac{\Gd}{\Gd J(x_m)}Z_J \right|_{J=0}.
$$

\subsection{Wick's theorem and Feynman graphs}
Wick's theorem now states that, similarly to \eqref{wick},
$$
   \left. \frac{\Gd}{\Gd J(x_1)}\dots\frac{\Gd}{\Gd J(x_m)}e^{\frac12\<J,GJ\>}
    \right|_{J=0}=\sum G(x_{i_1},x_{i_2})\dots G(x_{i_{m-1}},x_{i_m}),
$$
where the sum is over all pairings $(i_1,i_2)\dots(i_{m-1},i_m)$
of $1,\dots,m$. Just as in the finite dimensional case, we may
encode each pairing by a graph with $m$ univalent vertices
labelled by $1,\dots,m$, and edges connecting vertices $i_1$ with
$i_2$, \dots, and $i_{m-1}$ with $i_m$ presenting the factors of
$G$.

Let us add a potential \eqref{U} to the action and define
$$
    Z_U=\int \cD\phi\ e^{-1/2\<K\phi,\phi\>+\h U(\phi)\>}.
$$
Then, similarly to \eqref{Z_U}, we have
$$
        Z_U=Z_0 e^{\h U(\frac{\Gd}{\Gd J})}e^{\frac12\<J,GJ\>}\big|_{J=0}.
$$
Using again the Wick's theorem, we can rewrite the latter
expression in terms of Feynman graphs to get
\begin{equation}\label{graphs}
Z_U=\sum_{\GG}\frac{\h^{|\GG|}}{\Aut\GG}\int_{labels}\prod_vU_v\prod_e G_e ,
\end{equation}
where the integral is over all labellings of the ends of edges,
$U_v=U(x_1,\dots,x_k)$ for a $k$-valent vertex with the labels
$x_1,\dots,x_k$ of the adjacent edges, and $G_e=G(x_i,x_j)$ for an
edge with labels $x_i,x_j$. Sometimes it is convenient to include
the integration over the labels of the star of a vertex into the
weight of this vertex.

\subsection{An example: $\phi^4$-theory}
\begin{figure}
\includegraphics[height=0.75in]{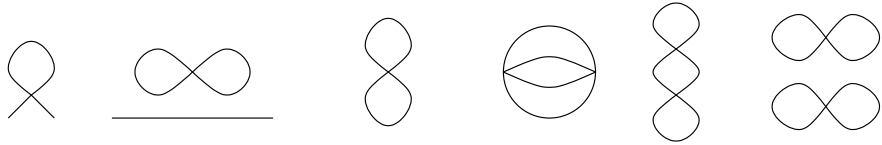}
\caption{Graphs appearing in the $\phi^4$-theory} \label{fig_phi4}
\end{figure}

Let us write down the Feynman rules for a potential $U(\phi)=\int
dx\ \phi^4(x)$. Firstly, the relevant graphs have vertices of
valence one or four. Secondly, all edges adjacent to a vertex
should be labelled by the same $x$, so we may instead label the
vertices. An edge with labels $x,y$ represents $G(x,y)$ and
(including the integration over the vertex labels into the weights
of vertices) an $x$-labelled vertex represents $\int dx$. The
linear term in the power series expansion of $\<x_1,x_2\>_U$
should correspond to non-vacuum graphs with two legs, labelled by
$x_1$ and $x_2$, and one 4-valent vertex. There is only one such
graph, see Figure \ref{fig_phi4}.
It represents $\int dx\ G(x_1,x)G(x,x)G(x,x_2)$ and should enter
with the multiplicity 12 (the number of all pairings of 6 vertices
$x_1,x_2,x,x,x,x$ in which $x_1$ is not connected to $x_2$). Let
us now check this directly. Indeed, the coefficient of $\h$ in
$\<x_1,x_2\>$ is
$$\frac{\Gd}{\Gd J(x_1)}\frac{\Gd}{\Gd J(x_2)}
\frac{\Gd^4}{\Gd J(x)^4} e^{\frac12\<J,GJ\>}\big|_{J=0}= 12\int
dx\ G(x_1,x)G(x,x)G(x,x_2),$$ where we applied Wick's theorem to
obtain the desired equality.

In a similar way, the linear term in the expansion of $Z_U$ should
correspond to the graph with no legs and one vertex of valence
four (see Figure \ref{fig_phi4}), representing $\int dx\ G^2(x,x)$
(and entering with the multiplicity 3).

\begin{exs}[$\phi^3$-theory]
Let $U(\phi)=\int dx\ \phi^3(x)$. Find the Feynman rules for this
theory. Which graphs will contribute to the coefficient of $\h^2$
in the power series expansion of the 2-point function
$\<x_1,x_2\>_U$? Write down these coefficients explicitly.
\end{exs}

\subsection{Convergence}
\label{sub:divergence}
Usually the integrals which we get by a perturbative Feynman
expansion are divergent and ill-defined in many ways. Often
one has to {\em renormalize} (i.e. to find some way to remove
divergencies in a unified and consistent manner) the theory to
improve its behavior. Until recently renormalization was considered
by mathematicians more like a physical art than a technique;
lately Connes and Kreimer \cite{CK} have done some serious work
to explain renormalization in purely mathematical terms
(see a paper by Kreimer in this volume).

But even in the best cases, the Green function $G(x,y)$ usually
blows up near the diagonal $x=y$, which brings two problems:
Firstly, the weights of graphs with looped edges, starting and
ending at the same point (so-called {\em tadpoles}) are ill-defined
and one has to get rid of them in one or another way. Secondly,
all diagonals have to be cut out from the spaces over which the
integration is performed, so the resulting configuration spaces
are open and the convergence of all integrals defining the weights
has to be proved. Mathematically these convergence questions
usually boil down to the existence of a Fulton-MacPherson-type (see
\cite{FMcP}) compactification of configuration spaces, to which
the integrand extends.

There is also a challenging problem to interpret the Feynman
diagrams series in some classical mathematical terms and to
understand the way to produce them without a detour to physics
and back. In many examples this may be done in terms of a
homology theory of some grand configuration spaces glued from
configuration spaces of different graphs along common boundary
strata.

We will see all this on an example of the Chern-Simons theory
in the next section.


\section{An example of QFT: Chern-Simons theory}
\label{sec:CS}

The Chern-Simons theory has an almost topological character and as
such presents an interesting object for low-dimensional topologists.
For a connection $\Ga$ in a trivial $SU(N)$-bundle over a 3-manifold
$M$ one may define (\cite{CS}; see also \cite{F}) the Chern-Simons
invariant $CS(\Ga)$ as described in Section \ref{sub:CStheory} below.
It is the action functional of the classical Chern-Simons theory and
was extensively used in mathematics for many years to study properties
of 3-manifolds (mostly due to the fact that the classical solutions,
i.e. the critical points of $CS(\Ga)$, are flat connections). But it
is the corresponding quantum theory which is of interest for us.
Its mathematical treatment started only about a decade ago, following
Witten's suggestion \cite{W} that it leads to some interesting invariants
of links and $3$-manifolds, in particular, to the Jones polynomial.
While Witten's idea was based on the validity of the path integral
formulation of the quantum Chern-Simons theory, his work catalyzed
much mathematical activity. By now mathematicians more or less managed
to formalize the relevant perturbative series and exorcize from them
all physical spirit, leaving a (surprisingly rich) rigorous mathematical
extract. In this section I will describe this process in a number of
iterations, starting from an intuitive and roughest description
and slowly increasing the level of rigor and details. Finally, I
will try to reinterpret these Feynman series in some classical
topological terms and formulate some corollaries.

\subsection{Chern-Simons theory}\label{sub:CStheory}
Further we will use the following data:
\begin{itemize}
\item A closed orientable 3-manifold $M$ with an oriented framed
      link $L$ in $M$.
\item A compact connected Lie group $G$ with an Ad-invariant trace
      $\Tr:\fG\to\R$ on the Lie algebra $\fG$ of $G$.
\item A principal $G$-bundle $\cP\to M$.
\end{itemize}

To simplify the situation, we will additionally assume that $G$ is
simply connected, since for such groups any principal $G$-bundle
over a manifold $M$ of dimension $\le 3$ (which is our case) is
trivializable, see e.g. \cite{F}.

The appropriate notions of the Chern-Simons theory, considered as
a field theory, are as follows. The manifold $M$ plays the role of
the space-time manifold $X$. Denote by $\cA$ the space of
$G$-connections on $\cP$ and let $\cG=\Aut(\cP)$ be the gauge
group. Fields $\Gp$ on $M$ are $G$-connections on $\cP$, i.e.
$\cF=\cA$.
The Lagrangian is a functional $L:\cA\to\GO^3(M)$ defined by
$$
   L(\Ga)=\Tr(\Ga\wedge d\Ga+\frac23\Ga\wedge\Ga\wedge\Ga).
$$
\begin{rem}

This choice can be motivated as follows.
Let $\Gt=d\Ga+\Ga\wedge\Ga$ be the curvature of $\Ga$.
Then $\Tr(\Gt\wedge\Gt)$ is the Chern-Weil 4-form\footnote{
Chern-Weil theory states that the de Rham cohomology class
of this form is a certain characteristic class of $\cP$}
on $\cP$, associated with $\Tr$; this form
is gauge invariant and closed. The Chern-Simons Lagrangian
$\displaystyle{CS(\Ga)=\Tr(\Ga\wedge\Gt+\frac23\Ga\wedge\Ga\wedge\Ga)}$
is an antiderivative of $\Tr(\Gt\wedge\Gt)$ on $\cP$:
it is a nice exercise to check that $d(CS(\Ga))=\Tr(\Gt\wedge\Gt)$.
\end{rem}

The corresponding Chern-Simons action is a function $CS:\cA\to\R$
given by
$$
CS(\Ga)=\frac{1}{4\pi}\int_M dx\Tr(\Ga\wedge d\Ga+\frac23\Ga\wedge\Ga\wedge\Ga).
$$
It is known that the critical points of this action correspond to flat
connections and (assuming that $\Tr$ satisfies a certain integrality
property\footnote{Namely that the closed form $\frac1{6\pi}
\Tr(\Ga\wedge\Ga\wedge\Ga)$ represents an integral class in
$H^3(G,\R)$}, which holds in particular for the trace in the fundamental
representation of $G$) it is gauge invariant modulo $2\pi\Z$.

The partition function is given by the following path integral:
\begin{equation}\label{eqZ}
Z=\int_{\cA} e^{ik CS(\Ga)}\cD\Ga.
\end{equation}
\no Here the constant $k\in\N$ is called {\em level} of the theory;
its integrality is needed for the gauge invariance of $Z$.

Now, let $L=\cup_{j=1}^m L_j$, $j=1,\dots m$ be an oriented framed
$m$-component link in $M$ such that each $L_j$ is equipped with a
representation $R_j$ of $G$. Given a connection $\Ga\in\cA$, let
$\hol_{L_j}(\Ga)$ be the holonomy
\begin{equation}\label{eq:hol}
      \hol_{L_j}(\Ga)=\exp\oint_{L_j} \Ga
\end{equation}
of $\Ga$ around $L_j$.
Observables in the Chern-Simons theory are so-called Wilson loops.
The {\em Wilson loop} associated with $L_j$ is the functional
$$
           \cW(L_j,R_j)=\Tr_{R_j}(\hol_{L_j}(\Ga)).
$$
The $m$-point correlation function $\<L\>=\<L_1,L_2,\dots L_m\>$
is defined by

\begin{equation}\label{eqW}
\<L\>=Z^{-1} \int_{\cA} e^{ikCS(\Ga)}
\prod_{j=1}^m{\cW}(L_j,R_j)\cD\Ga.
\end{equation}

Since the action is gauge invariant, extrema of the action
correspond to points on the moduli space of flat connections.
Near such a point the action has a quadratic term (arising
from $\Ga\wedge d\Ga$) and a cubic term (arising from
$\Ga\wedge\Ga\wedge\Ga$).
We would like to consider a perturbative expansion of this theory.

\subsection{What do we expect}
Which Feynman graphs do we expect to appear in the perturbative
Chern-Simons theory?

Firstly, a gauge-fixing has to be performed, so the ghosts have to
be introduced. As a result, we should have two types of edges:
the usual non-directed edges (corresponding to the inverse of the
quadratic part) and the directed ghost edges.

Secondly, in addition to the quadratic term the action contains a
cubic term, so the internal vertices should be trivalent.
Also, this time the cubic term is given by an antisymmetric tensor
instead of a symmetric one, so one should fix a cyclic order at each
trivalent vertex, with its reversal negating the weight of a graph.
Two types of edges should lead to two types of internal vertices:
usual vertices where three usual edges meet, and ghost vertices
where one usual edge meets one incoming and one outgoing ghost
edge.

Thirdly, note that the situation with legs is somewhat different
from our earlier considerations. Indeed, the legs (i.e. univalent
ends of usual edges) of Feynman graphs, instead of being fixed at
some points, should be  allowed to run over the link $L$,
with each link component entering in $\<L\>$ via its holonomy
\eqref{eq:hol}. To reduce this to our previous setting, we can
use Chen's iterated integrals to expand the holonomy in a power
series where each term is a polynomial in $\Ga$. In terms of a
parametrization $L_j:[0,1]\to\R^3$, this expansion can be written
explicitly using the pullback $L_j^*\Ga$ of $\Ga$ to $[0,1]$ via
$L_j$:

\begin{multline*}
\hol_{L_j}(\Ga)=1+\int_{0<t<1}(L_j^*\Ga)(t)+
\int_{0<t_1<t_2<1} (L_j^*\Ga)(t_2) \wedge (L_j^*\Ga)(t_1)+\\ \dots +
 \int_{0<t_1<\dots<t_k<1} (L_j^*\Ga)(t_k)\wedge \dots \wedge (L_j^*\Ga)(t_1)+\dots
\end{multline*}
where the products are understood in the universal enveloping
algebra $U(\fG)$ of $\fG$. Thus we should sum over all graphs with
any number $k_j$ of cyclically ordered legs on each $L_j$, and
integrate over the positions

 $$
    (x(t_1), \dots, x(t_{k_j}))\in  L_j^{k_j},\quad   0<t_1<\dots<t_{k_j}<1
$$ of these legs.

These simple considerations turn out to be quite correct. Of
course, one should still find an explicit formulas for the weights
of such graphs. An explicit deduction of the Feynman rules for the
perturbative Chern-Simons theory is described in details in
\cite{BN,GMM}. Let me skip these lengthy calculations and
formulate only the final results. For simplicity I will consider
only an expansion around the trivial connection in $M=\R^3$.

\subsection{Feynman rules}
\begin{figure}[htb]
\includegraphics[width=5.0in]{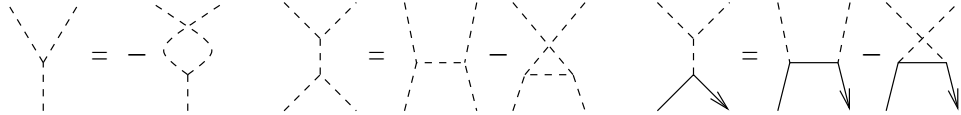}
\caption{Antisymmetry and Jacobi relations} \label{fig:rels}
\end{figure}

It turns out (see \cite{BN}) that the weight system
$W^{CS}$ of the perturbative Chern-Simons theory splits as
$W^{CS}=W^G W$, where $W^G$ contains all the relevant Lie-algebraic
data of the theory (but does not depend on the location of the
vertices of a graph), and $W$ contains only the space-time
integration. Since the whole construction should work for any Lie
algebra, one may encode the antisymmetry and Jacobi relations
already on the level of graphs, changing the weight $W_\GG^G$ of
a graph $\GG$ to a ``universal weight" $[\GG]$, which is an
equivalence class of $\GG$ in the vector space over $\Q$
generated by abstract (since we do not care about the location
in $\R^3$ of their vertices) graphs, modulo some simple
{\em diagrammatic antisymmetry and Jacobi relations}, shown on
Figure \ref{fig:rels}. The same relations hold for graphs with
either usual or ghost edges, so we may think that the relations
include the projection making all edges of one type.

The drawing conventions merit some explanation. It is assumed that
the graphs appearing in the same relation are identical outside
the shown fragment. In each trivalent vertex we fix a cyclic order
of edges meeting there; unless specified otherwise, it is assumed
to be counter-clockwise. The edges are shown by dashed lines, and
the link component $L_j$ (fixing the cyclic order of the legs) by
a solid line. An important consequence of the antisymmetry
relation is that for any graph $\GG$ with a tadpole (a looped edge)
we have $[\GG]=0$ due to the existence of a ``handle twisting"
automorphism, rotating the looped edge. Thus from the beginning we
can restrict the class of graphs to graphs without tadpoles.

It remains to describe the weight $W_\GG$ of a graph $\GG$.
Roughly speaking, for each internal vertex of $\GG$ we are to
perform integration over its position in $\R^3$ (for a ghost
vertex we should also take a certain derivative acting on the term
corresponding to the outgoing ghost edge); for each leg we are to
perform integration over its position in $L_j$ (respecting the
cyclic order of legs on the same component). As for the edges, we
are to assign to each usual and ghost edge inverses of the
operator $curl$ and of the Laplacian, respectively.

Somewhat surprisingly (see e.g. \cite{GMM}) two types of edges may
be neatly joined into one ``combined" edge, thus reducing the
graphs in question to graphs with only one type of edges (and just
one type of uni- and trivalent vertices). The weight
$G(x(e),y(e))$ of such an edge $e$ with the ends in $(x(e),y(e))
\in \R^3\times\R^3$ has a nice geometrical meaning: it is given by
$G(x,y)=\Go(x-y)$, where $$\Go(x)=\frac{x^1 dx^2 \wedge dx^3}
{2\pi||x||^3}+\text{cyclic permutations of (1,2,3)}$$ is the
uniformly distributed area form on the unit 2-sphere $|x|=1$ in
the standard coordinates in $R^3$. In fact, the usual and the
ghost edges (with two possible orientations) give respectively the
$(1,1)$, $(2,0)$, and $(0,2)$ parts of $\Go(x-y)$ in terms of its
dependence on $dx$ and $dy$. Abusing notation, I will depict the
combined edge again by a dashed line.

\begin{rem}
A simple explanation for an existence of such a simple unified
propagator escapes me. The only explanation which I know is way
too complicated: it is the existence (see \cite{AS}) of the
``superformulation" of the gauge-fixed theory, i.e. the fact
that the connection together with the ghosts may be united in
a ``superconnection" of a supertheory, which leads to an existence
of a ``superpropagator", uniting the usual and the ghost propagator.
I believe that there is a simple explanation, probably emanating
from the scaling properties and the topological invariance of the
Chern-Simons theory, by which one should be able to predict that
the combined propagator should be dilatation- and rotation-invariant.
\end{rem}

\begin{rem}
Note that the weight $G(x,x)$ of a tadpole is not well-defined,
so it is quite fortunate that we got $[\GG]=0$ for any such graph.
\end{rem}

To sum it up, we are interested in the value
\begin{equation}\label{eq:L}
\<L\>=\sum_\GG\frac{W_\GG\h^{|\GG|}}{\Aut\GG}[\GG]
\end{equation}
where $|\GG|$ is half of the total number of vertices (univalent
and trivalent) of $\GG$,
and the weight of $\GG$ is
given by the integral
\begin{equation}\label{eq:W}
W_\GG=\int_{C_\GG}\prod_e G(x(e),y(e))
\end{equation}
over the space $C_\GG$ of all possible positions of vertices of
$\GG$, such that all vertices remain distinct.
Here $\h=(k+h^\vee)^{-1}$, where $h^\vee$ is the dual Coxeter number of $G$
(see \cite{W}).

I shall describe in more details the type of graphs which appear
in this formula and their weights $W_\GG$ (both the
configuration spaces $C_\GG$, and the integrand).

\subsection{Jacobi graphs}
Let us start with the graphs. Instead of thinking about graphs
embedded in $\R^3$, consider abstract graphs (with just one type
of edges), such that
\begin{itemize}
\item all vertices have valence one (legs) or three;
\item there are no looped edges;
\item all legs are partitioned into $m$ subsets $l_1,\dots,l_m$;
\item legs of each subset $l_j$ are cyclically ordered;
\item each trivalent vertex is equipped with a cyclic order of
three half-edges meeting there;
\end{itemize}
for technical reasons it will be convenient to think that, in
addition to the above,
\begin{itemize}
\item all edges are ordered and directed.
\end{itemize}
We will further address the last three items simply as an
orientation of a graph.

\begin{figure}[htb]
\includegraphics[width=5.0in]{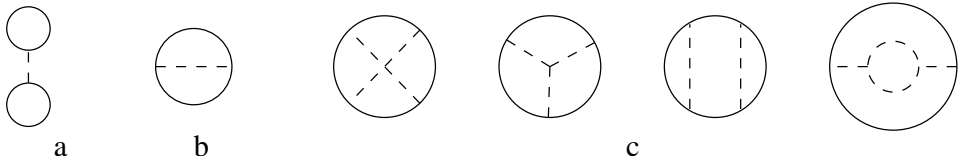}
\caption{Graphs of degree one and two} \label{fig:xyz}
\end{figure}

For such a graph $\GG$ with a total of $2n$ (univalent and
trivalent) vertices define the {\em degree} of $\GG$ by $|\GG|=n$,
and denote the set of all such graphs by $\bJ_n$. Set
$\bJ=\cup\bJ_n$. The ordering and directions of edges of
graphs in $\bJ$ may be dropped by an application of an obvious
forgetful map. See Figure \ref{fig:xyz} for graphs of
degree one with $m=2$ and $m=1$, and graphs of degree two with $m=1$.
Both antisymmetry and Jacobi relations of Figure \ref{fig:rels}
preserve the degree of a graph, thus we may consider a vector
space over $\Q$ generated by graphs in $\bJ_n$ modulo forgetful,
antisymmetry and Jacobi relations. We will call it the {\em space
of Jacobi graphs of degree $n$} and denote it by $\J_n$; denote
also $\J=\oplus_n\J_n$, and let as before $[\GG]$ be the class of
$\GG\in\bJ$ in $\J$.

\begin{exs}
Let $m=2$. Write the relations between the equivalence classes of
degree two graphs shown in Figure \ref{fig:xyz}c. What is the
dimension of $\J_2$?
\end{exs}

This settles the type of graphs appearing in formula \eqref{eq:L}:
the summation is over all graphs in $\bJ$, while
$\<L\>\in\J[[\h]]$. It is somewhat simpler to study separately the
components of different degrees; define
\begin{equation}\label{eq:Ln}
\<L\>_n=\sum_{\GG\in\bJ_n}\frac{W_\GG}{|\Aut(\GG)|}[\GG]
\end{equation}

\subsection{Configuration spaces}
Let us deal now with the weights \eqref{eq:W} of graphs
(see \cite{BT, Poi, T} for details).
The domain of integration in \eqref{eq:W} is the configuration space
$C_\GG$ of embeddings of the set of vertices of $\GG$ to $\R^3$,
such that the legs of each subset $l_j$ lie on the corresponding
component $L_j$ of the link $L$ in the correct cyclic order. It is
easy to see that for a graph $\GG$ with $k$ trivalent vertices and
$k_j$ legs ending on $L_j$, $j=1,\dots,m$ we have
$C_\GG\cong(\R^3)^k\times\prod_j (S^1 \times \Gs^{k_j-1})
\sminus\GD$, where $\Gs^k$ is a $k$-dimensional simplex, and $\GD$
is the union of all diagonals where two or more points coincide.
Indeed (forgetting for a moment about coincidences of vertices),
each trivalent vertex is free to run over $\R^3$, while $k_j$ legs
ending on $L_j$ run over $S^1\times\Gs^{k-1}$, where $S^1$ encodes
the position of the first leg, and the following legs are encoded
by their distance from the previous one.
\begin{exs}\label{exs:count}
Show that the dimension of $C_\GG$ is twice the number of the
edges of $\GG$.
\end{exs}

Now, an orientation of a graph $\GG$ determines an orientation of
$C_\GG$; its idea is in fact based on Exercise \ref{exs:count}.
Let me describe this construction in some local coordinates. Near
each trivalent vertex of $\GG$ there are three local coordinates
(describing its movement in $\R^3$); assign one of them to each of
the three ends edges meeting in this vertex using their cyclic
order. Near each leg of $\GG$ there is only one local coordinate
(describing its movement along the link); assign it to the
corresponding end of the edge. By now the end of any edge has one
coordinate assigned to it. It remains to order them using the
given ordering of all edges of $\GG$ and their directions. Let us
order them as $(x_1,y_1,x_2,y_2,\dots,x_n,y_n)$ where $(x_i,y_i)$
are the coordinates assigned to the beginning and the end of
$i$-th edge. This defines an orientation of $C_\GG$.

\begin{exs}
The above construction involves a choice in each trivalent vertex
since we had only a {\em cyclic} order of the edges meeting there,
while we used a total order of these three edges. Show that a
cyclic permutation of the three local coordinates used there
preserves the orientation of $C_\GG$. Also, we used the
orientation of $\GG$; what happens with the orientation of $C_\GG$
if:
\begin{enumerate}
\item
The cyclic order of three half-edges in one vertex is reversed?
\item
A pair of edges is transposed in the total ordering of all edges?
\item
The direction of an edge is reversed?
\end{enumerate}
\end{exs}

\subsection{Gauss-type maps of configuration spaces}
To understand the integrand in \eqref{eq:W}, consider a directed
edge $e$. Its ends $(x,y)$ represent a point in the square
$\R^3\times\R^3$ with the diagonal $\GD=\{(x,y)|\ x=y\}$ cut out.
This cut square $C=\R^3\times\R^3\sminus\GD$ has the homotopy type
of $S^2$, with the Gauss map $$\phi:(x,y)\mapsto \frac{y-x}
{||y-x||}$$ providing the equivalence. The form $\Go(y-x)$
assigned to this edge is nothing more than a pullback of the area
form $\Go$ on $S^2$ to $C$ via the Gauss map: $$\Go(y-x)=\phi^*
\Go$$ Each edge $e$ of a graph $\GG$ defines an evaluation map
$\ev_e:C_\GG\to C$, by erasing all vertices of $\GG$ but for the
ends of $e$. The composition $\phi_e=\phi\circ\ev_e$ defines the
Gauss map corresponding to (the ends of) an edge $e$. The graph
$\GG$ with an ordering $e_1,\dots,e_n$ of edges defines the
product
$\phi_\GG=\prod_{i=1}^n\phi_{e_i}:C_\GG\to(S^2)^n$ of Gauss maps.
Finally, the weight $W_\GG$ is given by integrating the pullback
of the volume form $\dvol =\wedge_{i=1}^n\Go$ on $(S^2)^n$ to
$C_\GG$ by the product Gauss map $\phi_\GG$:
\begin{equation}\label{eq:degree}
W_\GG=\int_{C_\GG}\phi_\GG^*\, \dvol
\end{equation}

\begin{exs}\label{exs:double}
Suppose that a graph $\GG$ has a double edge (i.e., a pair of
edges both endpoints of which coincide). Show that
$\dim(\phi_\GG(C_\GG))\le\dim(C_\GG)-1$. Deduce that $W_\GG=0$.
\end{exs}

The following important example shows that at least in some simple
cases $W_\GG$ has an interesting topological meaning:

\begin{ex}\label{ex:link}
Let $\GG=e$ be a graph with one edge with the ends on two link
components $L$, see Figure \ref{fig:xyz}a. The configuration
space $C_e\cong S^1\times S^1\subset C$ is a torus. It is mapped
to $S^2$ by the Gauss map $\phi=\phi_e$. The weight $W_e=\int_{C_e}
\phi^*\Go=\deg(\phi)$ is in this case just the degree of the map
$\phi$. This fact has many important consequences. In particular
$W_e$ takes only integer values and is preserved if we change the
uniformly distributed area form $\Go$ to any other volume form
$\dvol$ on $S^2$ normalized by $\int_{S^2}\dvol=1$. It is also
preserved if we deform the link by isotopy (since then the
configuration space changes smoothly and the degree can not jump),
so is a link invariant. This invariant is easy to identify:
$\int_{C_e}\phi^*\Go=\lk(L_1,L_2)$ is the famous Gauss integral
formula for the linking number $\lk(L_1,L_2)$ of $L_1$ with $L_2$.
Thus we get
\begin{prop}
Let $\GG=e$ be a graph with one edge with the ends on two link
components $L$. Then the weight $W_e$ is the linking number
$\lk(L_1,L_2)$.
\end{prop}
\end{ex}

\begin{figure}[htb]
\includegraphics[height=0.75in]{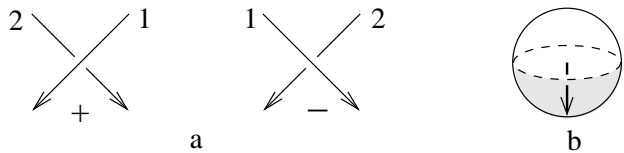}
\caption{Signs of crossings and the south pole on $S^2$}
\label{fig:signs}
\end{figure}

\begin{exs}
There is a simple combinatorial way to compute $\lk(L_1,L_2)$
from any link diagram: count all crossings where $L_1$ passes
over $L_2$, with signs shown in Figure \ref{fig:signs}a.
Interpret this formula as a calculation of $\deg(\phi)$ by
counting (with signs) the number of preimages of a certain
regular value of $\phi$ (hint: look at Figure~\ref{fig:signs}b).
What formula would we get if we  counted the preimages
of the north pole?
\end{exs}

For other graphs the situation is more complicated. For example,
let $\GG=e$ be the graph with one edge, both ends of which end on
{\em the same} link component, see Figure~\ref{fig:xyz}b. Then the
configuration space $C_e$ is an {\em open} annulus $(R^3)^0 \times
S^1\times\Gs^1\smallsetminus\GD=S^1\times(0,1)$ (torus cut along
the diagonal). The Gauss map $\phi$ is badly behaved near the
diagonal, so the integrand blows up near the diagonal and we can
not extend it to the closed torus. The integral nevertheless
converges; one way to see it is to compactify $C_e$, cutting out of
it some small neighborhood of the diagonal. This makes $C_e$ into
a closed annulus $C^\Ge_e=S^1\times [\Ge,1-\Ge]$ (thus making the
integral convergent) and we can recover the initial integral by
taking $\Ge\to 0$. But the Gauss integral $W_e$ is no more
a knot invariant: it may take any real value under a knot isotopy.
A detailed discussion on this subject may be found in \cite{BT}.
Why does this happen? The reason is that the compactified space
$C_e$ is not a torus, but an annulus, so has a boundary and the
degree of the Gauss map is not well-defined. When both ends of
the edge start to collide together, the direction of the vector
connecting them (which appears in the Gauss map) tends to the
(positive or negative) tangent direction to the knot. The image
of the unit tangent to the knot under the Gauss map is a certain
curve $\Gg$ on $S^2$.
One of the boundary circles $S^1\times\Ge$ and $S^1\times(1-\Ge)$
of $C^\Ge_e$ is mapped into $\Gg$, while the other is mapped into
$-\Gg$, and the weight $W_e$ is part of the area of $S^2$ covered
by the annulus $\phi(C^\Ge_e)$ between these curves. Unfortunately,
$\Gg$ may move on $S^2$ under an isotopy of $L$, so this area may
change.

In this particular case there is a neat way to solve this problem:
let $L$ be {\em framed} (i.e. fix a section of its normal bundle).
We may think about the framing as about a unit normal vector
$n(x)$ in each point $x$ of a knot. This allows us to slightly
deform the Gauss map: $\phi_(x,y)\to\phi(x,y)+\Ge n(y)$. Now
both boundary circles of the annulus $C_e$ map into the same curve
on $S^2$ (why?) and we may glue the annulus into the torus so that
the map $\phi_\Ge$ extends to it. It makes $W_e$ into an invariant
of framed knots, called the self-linking number (the same result
may be obtained by slightly pushing $L$ off itself along the
framing and considering the linking number of the knot with its
pushed-off copy).

It turns out that for other graphs there are also no divergence
problems, so all integrals $W_\GG$ converge, and that a collision
of all vertices of a graph to one point (so-called anomaly, see
\cite{Poi, T}) is the only source of non-invariance, exactly as
for $W_e$ above.
Thus there is a suitable normalization of the expression \eqref{eq:Ln}
for $\<L\>=\sum_n\<L\>_n$ which gives a link invariant.
To avoid a complicated explicit description of this normalization,
let me formulate this result as follows:

\begin{thm}[\cite{AF}, \cite{Poi}, \cite{T}]\label{thm:L}
Let $L=\cup_{i=1}^m L_i$ be a link. Then $\<L\>$ depends only on
the isotopy class of $L$ and on the Gauss integrals $W_e(L_i)$
of each component $L_i$. In particular, an evaluation of $\<L\>$
at representatives of $L$ for which $W_e(L_1)=\dots=W_e(L_m)=0$
is a link invariant.
\end{thm}

\begin{rem}
It is known that this is a universal invariant of finite type.
In particular this means that it is stronger than both the Alexander
and the Jones polynomials (it contains the two-variable HOMFLY
polynomial) and all other quantum invariants. Conjecturally the
anomaly vanishes and this invariant coincides with the Kontsevich
integral, see \cite{Poi}.
\end{rem}

\begin{ex}\label{ex:v2}
Let $L$ be a knot, and take $n=2$. There are four graphs of degree
two, shown in Figure \ref{fig:xyz}c. We will denote the first of
them $X$, and the second by $Y$. By Exercise \ref{exs:double} the
weight of the third graph vanishes. Also, choose a framing of $L$
so that the self-linking is 0; then the contribution of the last
graph vanishes (another way to achieve the same result is to add
to $\<L\>_2$ a certain multiple of the self-linking number
squared); we can set then $[X]=[Y]$. Thus we will consider simply
$$
    v_2=\frac14\int_{C_X} \phi^*_X (\Go\wedge\Go)+\frac13\int_{C_Y}\phi^*_Y
(\Go\wedge\Go\wedge\Go).
$$
The first integral is 4-dimensional,
while the second is 6-dimensional; none of them separately is a
knot invariant (see \cite{PV} for a discussion); however, their
sum $v_2$ is (see \cite{BN}, \cite{PV})! This invariant is, up
to a constant, the second coefficient of the Alexander-Conway
polynomial. See \cite{PV} for its detailed treatment as the
degree of a Gauss-type map.
\end{ex}

\subsection{Degrees of maps}
How can we explain the result of Theorem \ref{thm:L}? We may try
to repeat the reasoning of Example \ref{ex:link} in the general
case. Recall that by Exercise \ref{exs:count} the dimensions of
$C_\GG$ and $(S^2)^n$ match, so if $C_\GG$ would be a closed
manifold, then equation \eqref{eq:degree} would be a formula for a
calculation of the degree of $\phi_\GG$. In other words, if
$C_\GG$ would have a fundamental class, \eqref{eq:degree} would be
its pairing with the pullback $\phi_\GG^*\dvol$. That would be
great: we would know that $W_\GG$ takes only integer values, and
would be able to compute it in many ways, including a simple
counting of preimages of any regular value of $\phi$.

Unfortunately, the reality is much worse: $C_\GG$ is an open
space, so the degree is not well-defined and even the convergence
is unclear. To guarantee the convergence, we should construct a
compactification $\bC_\GG$ of $C_\GG$ to which $\phi_\GG$ extends.
This however will cause new problems: the space $\bC_\GG$ will
have many boundary strata. There are various way to deal
with them: we can relativize some of them (i.e. consider a
relative version of the theory), cap-off some others (i.e.
glue to them some new auxiliary configuration spaces), or
zip them up (gluing a stratum to itself by an involution).
But in general, some boundary strata will remain; indeed, in Example
\ref{ex:v2} we have seen that none of $W_X$ or $W_Y$ separately
can be made into a knot invariant. The remedy would be to glue
together the configuration spaces for different graphs in
$\bJ_n$ along the common boundary strata. This tedious work
can be done indeed \cite{Poi, T} and (up to a certain anomaly
correction) one can interpret $\<L\>$ as the degree of a certain
map $\Phi_n$ from a grand configuration space $\cC_n$ to $(S^2)^n$.
Too many technicalities are involved to describe this construction
in necessary details, so I refer the interested reader to
\cite{Poi, T} and will present only a brief sketch of this
construction.

The first problem is that initially the dimensions of $C_\GG$ for
various $\GG\in\bJ_n$ do not match. E.g., in Example \ref{ex:v2},
for $n=2$ the spaces $C_X$ and $C_Y$ have dimensions 4 and 6
respectively.
This can be fixed by considering a product $C_\GG\times(S^2)^k$
of $C_\GG$ with enough spheres to make the maps $\phi_\GG\times
(\text{id})^k$ to have the same target space $(S^2)^N$ for all $\GG$.
Now one should do the gluings. When two endpoints of an edge $e$
of a graph $\GG$ collide, the corresponding boundary stratum of
$C_\GG$ looks like $C_G\times S^2$ for $G=\GG/e$.
Thus we can glue together such strata for all pairs $(\GG,e)$
with isomorphic $G=\GG/e$. Some more, so-called hidden, strata
remain after these main gluings. Fortunately, each of them
can be zipped-up (i.e. glued to itself by a certain involution).
The only codimension one boundary strata which remains
after all these gluings are the anomaly strata, where all vertices
of a graph $\GG\in\bJ$ collide together. These problematic anomaly
strata can be glued \cite{Poi} to a new auxiliary space. One ends
up with a grand configuration space $\cC_n$ endowed with a map
$\Phi:\cC_n\to (S^2)^N$ (glued from the corresponding maps
$\phi_\GG : \bC_\GG\rightarrow (S^2)^N$).
One may show that the cohomology $H^{2N}(\cC_n)$ of this space
projects surjectively to $\J_n$ (see \cite{KT} for a similar
case of 3-manifold invariants).
Then $\<L\>^0_n=\<L\>_n+\text{anomaly correction}$ can be
interpreted as the degree of $\Phi_n$, or more exactly, the
image in $\J_n$ of the fundamental class $[(S^2)^N]$ under
the induced composite map
$\pi\circ\Phi^*:H^{2N}(S^2)^N\to H^{2N}(\cC_n)\to\J_n$.

\subsection{Final remarks}
There remain many questions: which compactification should we
take, why do the antisymmetry and Jacobi relations appear in
the cohomology of the grand configuration space, etc.
Each of them is quite lengthy and is out of the scope of this
note. We refer the interested reader to \cite{BT, Poi, T}.
A mathematical treatment of invariants of 3-manifolds arising
from the Chern-Simons theory was done in \cite{AS, BC}; I
especially recommend \cite{KT}. While I do not know whether
similar Feynman series arising in other topological problems
always have a reformulation in terms of degrees of maps of
some grand configuration space, it seems quite plausible.
There are at least some other notable examples, see e.g.
\cite{P} for a similar interpretation of Kontsevich's
quantization of Poisson structures.

\end{document}